\documentclass[11pt]{article}
\usepackage{amsthm,amssymb,amsmath,amscd}
\usepackage[all]{xy}
\usepackage{mathrsfs}
\usepackage{multirow}

\usepackage{mathptm}

\newtheorem{Def}{Definition}[section]

\newtheorem{Theo}[Def]{Theorem}
\newtheorem{Lem}[Def]{Lemma}
\newtheorem{Koro}[Def]{Corollary}

\newtheorem{Rem}[Def]{Remark}

\newcommand{\add}{{\rm add}}

\newcommand{\Hom}{{\rm Hom }}
\newcommand{\pd}{{\rm pd }}%\newcommand{\pd}{{\rm proj.dim}}
\newcommand{\id}{{\rm injdim }}

\newcommand{\rd}{{\rm repdim }}
\newcommand{\gd}{{\rm gldim}}
\newcommand{\fd}{{\rm findim}}
\newcommand{\End}{{\rm End}}

\newcommand{\modcat}[1]{#1\mbox{{\rm -mod}}}
\newcommand{\Modcat}[1]{#1\mbox{{\rm -Mod}}}

\newcommand{\lra}{\longrightarrow}
\newcommand{\ra}{\rightarrow}

\begin{document}

{\Large \bf
\begin{center}
Representation dimensions linked by Frobenius bimodules
with applications to group algebras
\end{center}}
\medskip

\centerline{{\bf Changchang Xi}}

\begin{abstract}
We establish relations between representation dimensions of two algebras connected by a Frobenius bimodule or extension. Consequently, upper bounds and equality formulas for representation dimensions of group algebras, symmetric separably equivalent algebras and crossed products are obtained. Particularly, for any subgroup $H$ of a finite group $G$, if $[G:H]$ is invertible in a field, then the representation dimensions of the group algebras of $G$ and $H$ over the field are the same.
\end{abstract}

\renewcommand{\thefootnote}{\alph{footnote}}
\setcounter{footnote}{-1} \footnote{2020 Mathematics Subject
Classification: 16E20, 20C05, 16T05;
 18G20, 16S34.}
\renewcommand{\thefootnote}{\alph{footnote}}
\setcounter{footnote}{-1} \footnote{Keywords: Frobenius bimodule; Group algebra; Representation dimension; Separable equivalence; Sylow p-subgroup; Symmetric separable equivalence.}

\section{Introduction}
Let $A$ be an Artin algebra, that is, a $k$-algebra over a commutative Artin ring $k$ such that it is a finitely generated $k$-module. We denote by $A\modcat$ the category of all finitely generated left $A$-modules. For an  $X\in A\modcat$, we denote by $\add(X)$ the additive category generated by $M$ in $A\modcat$. An $A$-module $_AM\in A\modcat$ is called a \emph{generator} if ${}_AA \in \add(M)$; and \emph{cogenerator} if $D(A_A)\in \add(M)$, where $D$ stands for the usual duality for Artin algebras. Now, we first review Auslander's representation dimension of algebras in \cite{auslander}.

Let $A$ be an Artin algebra. The
representation dimension of $A$, denoted by $\rd(A)$, is defined to be
$$\rd(A):=\emph{inf}\, \{\gd(\End_A(M))\mid M={}_AA\oplus D(A_A)\oplus X,\, X\in A\modcat\, \}$$
where $\gd$ stands for the global dimension.

Clearly, $\rd(A)=0$ if and only if $A$ is semisimple, and there are no algebras $A$ with $\rd(A)=1$. Auslander proved in \cite{auslander} that $\rd(A)\le 2$ if and only if $A$ is representation-finite. It was shown in \cite{iyama} that $\rd(A)$ is finite for any Artin algebra $A$.

In this note, we establish relations between representation dimensions of algebras linked by symmetric separable equivalences or Frobenius bimodules, including Frobenius extensions.

Frobenius extensions were initially introduced by Kasch in \cite{kasch54}, and developed later by many authors (see the references in \cite{kadison}). As is known, Frobenius bimodules and extensions have many applications in different branches in mathematics, for instance, in knot theory, Yang-Baxter equations, code theory and representation theory of algebras (see \cite{g-wood, kadison, x1}).

In general, there are no expected relations between representation dimensions of algebras linked by arbitrary Frobenius bimodules. This can be seen from Frobenius extensions. On the one hand, if $A$ is a Frobenius $k$-algebra over a field $k$, then $k\hookrightarrow A$ is a Frobenius extension. In this case, $0=\rd(k)\le \rd(A)$, but the latter can be arbitrary. On the other hand, the $n\times n$ matrix algebra $M_n(k)$ over $k$ is a separable, Frobenius extension of the centrosymmetric matrix algebra $S_n(k)$ (see \cite{xy}). In this case, if $k$ has  characteristic $2$, then $\rd\big(S_n(k)\big)=2> 0=\rd\big(M_n(k)\big)$. Nevertheless, we shall establish the following two results in this direction. The first one is under some conditions on Frobenius bimodules.

\begin{Theo} \label{frob-bim} Let $A$ and $B$ be Artin algebras, and let $_BM_A$ be a Frobenius bimodule such that $_BM$ is a generator for $B\modcat$, $B$ is $M$-separable over $A$, and the $A$-$A$-bimodule $\End_B(M)$ is centrally projective over $A$. If $A$ is non-semisimple, then $\rd(B)\le \rd(A)$.
\end{Theo}

Before we state the second result related to representation dimensions and Frobenius bimodules, we first recall some definitions (see \cite{kadison93, kadison95, be, linckelmann, peacock}).

Let $R$ and $S$ be rings with identity. We say that $R$ separably divides $S$ if there are bimodules $_RP_S$ and $_SQ_R$ which are finitely generated and projective as one-sided modules, and there is a split epimorphism $\nu: {}_RP\otimes_S Q_R\ra {}_RR_R$ of $R$-$R$-bimodules \cite{kadison95}; and that $R$ symmetric separably divides $S$ if $R$ separably divides $S$, and
additionally, the functors $F:= {}_RP\otimes_S -$ and and $G:= {}_SQ \otimes_R-$ are an adjoint pair $(F,G)$ between the categories of all left modules, $S$-Mod and $R$-Mod, with $\nu$ the counit of adjunction.

Recall that rings $R$ and $S$ are \emph{separably equivalent} if $R$ separably divides $S$ with
split epimorphism  $\nu : {}_RP\otimes_S Q_R \ra {}_RR_R$ of $R$-$R$-bimodules and $S$ separably divides $R$ with split epimorphism  $\mu : {}_SQ\otimes_R P_S\ra {}_SS_S$ of $S$-$S$-bimodules. Further, rings $R$ and $S$ are \emph{symmetric separably equivalent} if $R$ symmetric separably
divides $S$ with adjoint functors $(F:= {}_RP\otimes_S-, G:= {}_SQ\otimes_R -)$ and $S$ symmetric separably
divides $R$ with adjoint functors ($G, F$) (see \cite[Definition 6.1, Remark 6.1]{kadison95}). In other words, the functors $F$ and $G$ between $S\Modcat$ and $R\Modcat$ are Frobenius functors \cite{morita}, that is, adjoint in either order. In this case, $_RP_S$ and $_BQ_A$ are Frobenius bimodules.

An example of symmetric separable equivalences is a separable, Frobenius extension $B\subseteq A$ of Artin algebras with a Frobenius homomorphism $E: A\ra B$ such that the restriction of $E$ to $B$ is the identity map $id_B$. Then $A$ and $B$ are symmetric separably equivalent by \cite[p.348]{Kadison2019} or \cite[Proposition 6.1]{kadison95}. Another example is a separable equivalence between symmetric Artin algebras $A$ and $B$. Recall that an Artin algebra $A$ is \emph{symmetric} if $A\simeq D(A)$ as $A$-$A$-bimodules where $D$ is the usual duality for Artin algebras. In this case, $A$ and $B$ are, in fact, symmetric separably equivalent by \cite[Corollary 5.5]{Kadison2019} for rings, or by \cite[Section 2]{peacock} for finite-dimensional algebras

\begin{Theo} Suppose that Artin algebras $A$ and $B$ are symmetric separably equivalent defined by bimodules $_AP_B$ and $_BQ_A$. If $_AP\otimes_BQ_A$ and $_BQ\otimes_AP_B$ are centrally projective over $A$ and $B$, respectively, then $\rd(A)=\rd(B).$
\label{sym-equ}
\end{Theo}

Our results can be applied to Frobenius extensions, and we then get two types of conditions for  inequality of representation dimensions. For details, we refer the reader to Corollary \ref{frob-ext}.
As to crossed products and group algebras, we mention here the following corollary which seems to be of interest for calculations of representation dimensions for group algebras of finite groups.

\begin{Koro} \label{int-cor}
$(1)$ For an Artin algebra $k$, and for any subgroup $H$ of a finite group $G$ with $[G:H]$ invertible in $k$, we have $\rd(k[G])\le\rd(k[H])$, where
$k[G]$ stands for the group algebra of $G$ over $k$. Particularly, if $k$ is a field,
then $\rd(k[G])$ = $\rd(k[H])\le |H|$.

$(2)$ Let $k$ be a field, and let $H$ be a finite-dimensional cocommutative Hopf-algebra over $k$, $A$ a left $H$-module algebra with a central trace $1$ element, and $\sigma$ an invertible normal two cocycle. Then $\rd(A\#_{\sigma}H) \le \rd(A)$, where $A\#_{\sigma}H$ is the crossed product of $A$ with $H$ by $\sigma$.

%$(3)$
\end{Koro}

Thus, for the group algebra $k[G]$ of a finite group $G$ over a field $k$ of characteristic $p$, if $P$ is any Sylow $p$-subgroup of $G$, then  $\rd(k[G])=\rd(k[P])\le |P|.$

\medskip
Proofs of all results in this section will be implemented in the next section where terminology unexplained in this section will be recalled.

\section{Proofs of the results}

Let $R$ be a ring with identity. We denote by $R\Modcat$ (respectively, $R\modcat$) the category of all (respectively, finitely generated) left $A$-modules. For an $R$-module $M$ in $R\modcat$, let $\add(M)$ be the full subcategory of $R\modcat$ consisting of all direct summands of direct sums of finitely many copies of $M$.
Composition of two homomorphisms $f:X\ra Y$ and $g: Y\ra Z$ will be denoted by $fg: X\ra Z$.

Let $S$ be another ring with identity and $_SM_R$ an $S$-$R$-bimodule. For convenience, we write $^*M$ for $\Hom_S(_SM, {}_SS_S) $ and $M^*$ for $\Hom_{R^{op}}(_SM_R,{}_RR_R).$ Recall that ${}_SM_R$ is called a \emph{Frobenius bimodule} if $_SM$ and $M_R$ both are finitely generated projective modules and ${}^*M \simeq M^*$ as $R$-$S$-bimodules; and $S$ is said to be \emph{$M$-separable over $R$} (see \cite{sugano1}) if the evaluation map $M\otimes_R{}^*M\to S$ is  split surjective as $S$-$S$-bimodules, that is, $_SS_S\in \add(_SM\otimes_R{}^*M).$ Note that if $M$ is a left $S$-generator, then $S$ is $M$-separable over $\End_S(M)$ since $ _SM\otimes_{\End_S(M)}{}^*M\simeq {}_SS_S.$ According to Hirata \cite{hirata}, an $R$-$R$-bimodule $N$ is said to be \emph{centrally projective} over $R$ if $_RN_R\in \add(_RR_R).$

\begin{Def} Let $f:S\ra R$ be an injective homomorphism of rings. We say that $f$ is

$(1)$ a \emph{Frobenius extension} (see \cite{kasch1}) if $_SR_R$ (or $_RR_S$) is a Frobenius $S$-$R$-bimodule (or $R$-$S$-bimodule);

$(2)$ a \emph{semisimple extension} if the multiplication map $R\otimes_SX\ra X$ is split surjective for all $R$-module $X$, or equivalently, the kernel of the map is a direct summand of $R\otimes_SX$ as $R$-modules.

$(3)$ a \emph{separable extension} if the multiplication map $R\otimes_SR\ra R$ is split surjective as $R$-$R$-bimodules;

$(4)$ an \emph{$H$-separable extension} if the $R$-$R$-bimodule $R\otimes_SR$ is centrally projective over $R$, that is, $R\otimes_SR\in \add(_RR_R)$ (see \cite{hirata}), and

$(5)$ a \emph{split extension} if $f$ is a split injective homomorphism of $S$-$S$-bimodules, that is, $_SS_S\in \add(_SR_S)$.
\end{Def}

Clearly, separable extensions are semisimple extensions. It is shown in \cite{hirata, sugano3} that $H$-separable extensions are separable extensions.
Also, Frobenius extensions generalise the notion of Frobenius algebra over a field. Recall that a finite-dimensional $k$-algebra $A$ over a field $k$ is called a \emph{Frobenius algebra} if $_AA\simeq \Hom_k(A_A,k)$ as left $A$-modules, that is, the map $k\ra A$ is a Frobenius extension.

A characterization of Frobenius extensions reads as follows: An extension $S\subseteq R$ of rings is  Frobenius if and only if there is an $S$-$S$-bimodule homomorphism $E:R\ra S$, $2n$ elements $x_i,y_i\in R$, $1\le i\le n$, such that, for any $a\in R,$
$$\sum_{i=1}^n x_i (y_ia)E = a = \sum_i (ax_i)E\, y_i.$$ In this case, $(E, x_i,y_i)$ is called a \emph{Frobenius system}, $E$ is called a \emph{Frobenius homomorphism} and $\sum_ix_iy_i$ is called a \emph{relative Casimir element} of the extension.
For further information on Frobenius extensions, we refer to \cite{kasch1, kadison}, and pertaining to separable extensions and their generalizations, we refer to \cite{zhu}.

On relations between $H$-separable extensions and Frobenius extensions we cite results in \cite[Proposition 2]{sugano1} (see also \cite[Proposition 3.2]{cun}), \cite[Theorem 7]{sugano1} and \cite[Theorem 1.2]{sugano2}, and \cite[Theorem 2]{sugano2} as the following lemma for the convenience of the later use.

\begin{Lem}\label{sep-over} $(1)$ Let $_SM_R$ be an $S$-$R$-bimodule such
that $M_R$ is finitely generated and projective, and let $E:=\End_{R^{op}}(M_R)$. Then $S$ is
$M$-separable over $R$ if and only if there is an $S$-$S$-bimodule homomorphism $h: E\to S$ with $(id_M)h=1$, that is, $_SS_S\in \add(_SE_S).$

$(2)$ Let $S\subseteq R$ be a split, $H$-separable extension of rings.

(a) If $_SR$ is finitely generated and projective, then $R\subseteq \End(_SR)$, $a\mapsto (x\mapsto xa),$ is a separable extnsion.

(b) If both $_SR$ and $R_S$ are finitely generated and projective, then $S\subseteq R$ is a Frobenius extension.

$(3)$ If $S\subseteq R$ is a separable extension such that $_SR_S$ is centrally projective over $S$, then $S\subseteq R$ is a Frobenius extension.
\end{Lem}

A typical example satisfying the conditions of Lemma \ref{sep-over}(1) is the following case: Let $k$ be a commutative ring, $G$ a group and $H$ a subgroup of $G$ such that $[G:H]<\infty$. We write $k[G]$ for the group ring of $G$ over $k$. If we define $S:=k[H]$, $R=k[G]$ and $M = {}_Sk[G]_R$, then $S$ is $M$-separable over $R$ by Lemma \ref{sep-over} since $_S\End_R^{op}(M_R)_S = _{k[H]}k[G]_{k[H]}$ and the decomposition of $G$ into double cosets of $H$ shows that $k[H]$ is a direct summand of $k[G]$ as $k[H]$-bimodules. Thus $k[H]$ is $_{k[H]}k[G]_{k[G]}$-separable over $k[G],$ that is, the extension $S\subseteq R$ is $H$-separable.

We need the following result, due to Auslander \cite{auslander}, which is based on the fact that the additive categories $\add(_AM)$ and $\add(_{\End_A(M)}\End_A(M))$ are equivalent under the functor $\Hom_A(M,-)$.

\begin{Lem} \label{lem1} Let $A$ be an Artin algebra and $M$ a generator-cogenerator for $A\modcat$. Define $E:=\End_A(M)$. Suppose $m\ge 2$ is an integer. Then the following are equivalent:

$(1)$ $\gd(E)\le m.$

$(2)$ $\pd_{E}\big(\Hom_A(M,Y)\big)\le m-2\,$ for all $Y\in A\modcat$, where $\pd$ means projective dimension.

$(3)$ For each $A$-module $Y\in A\modcat$, there is an exact sequence
$$0\lra M_{m-2}\lra ... \lra M_1\lra M_0\lra Y\lra 0,  $$
in $A\modcat$ with $M_j \in$ add$(_AM)$ for $j=0, ...., m-2$, such that the induced sequence
$$0\lra \Hom_A(X,M_{m-2})\lra ...
\lra \Hom_A(X,M_1)\lra \Hom_A(X,M_0)\lra
\Hom_A(,Y)\lra 0$$ is exact for all $X\in$ add$(_AM),$ or equivalently, $$0\lra \Hom_A(M,M_{m-2})\lra ...
\lra \Hom_A(M,M_1)\lra \Hom_A(M,M_0)\lra
\Hom_A(M,Y)\lra 0$$ is exact.
\label{replem}
\end{Lem}

\medskip
{\bf Proof of Theorem \ref{frob-bim}.} We define ${}_AN_B:= {}^*M=\Hom_B(_BM,{}_BB)$. Then $$M\otimes_A- \simeq \Hom_A(_AN,-):  A\mbox{-mod} \lra B\mbox{-mod} \quad \mbox{ and }\quad N\otimes_B- \simeq \Hom_B(M,-): B\mbox{-mod} \lra A\mbox{-mod}$$ (for example, see \cite{x1}). To prove Theorem \ref{frob-bim}, we show the following statements.

(1) If $_AX$ is a generator for $A$-mod, then $_BM\otimes_AX$ is a generator for $B$-mod because $_BM$ is a generator for $B$-mod.

(2) If $_AX$ is a cogenerator for $A$-mod, then $_BM\otimes_AX$ is a cogenerator for $B$-mod. In fact, if we assume $X=D(A_A)\oplus X'$, then $_BM\otimes_AX={}_BM\otimes_AD(A_A)\oplus {}_BM\otimes_AX'$. So we only need to prove that $_BM\otimes_AD(A_A)$ is a cogenerator for $B$-mod. But this follows from
$$ _BM\otimes_AD(A_A)\simeq \Hom_A(_AN_B,D(A_A))\simeq D(A\otimes_AN_B)\simeq D(N_B)=D(\Hom_{B}(_BM,{}_BB_B))$$and the fact that ${}_BB\in\add(_BM)$ implies
$$B_B\in \add(\Hom_B(_BM,{}_BB_B)) \mbox{ and } D(B_B)\in \add(D\Hom_B(_BM,{}_BB_B)) = \add\big(_BM\otimes_AD(A_A)\big).$$Hence ${}_BM\otimes_AX$ is a cogenerator for $B$-mod.

(3) Let $_AX$ be a generator-cogenerator for $A$-mod such that $n:=\rd(A)=\gd(\End_A(X))\ge 2$. Then we show $\gd(\End_B(M\otimes_AX))\le n.$

In fact, since $n:=\rd(A)=\gd(\End_A(X))\ge 2$, this means by Lemma \ref{lem1} that,
for any $A$-module $U$, there is an exact sequence $$ 0\lra U_{n-2}\lra U_{n-3}\lra \cdots\lra U_1\lra U_0\lra U\lra 0$$ such that $U_j\in\add(_AX)$ and $\Hom_A(X,-)$ preserves exactness of this sequence. Now, let $_BV$ be a $B$-module. Then, by Lemma \ref{lem1}, we have an exact sequence of $A$-modules:
$$ (*)\qquad 0\lra X_{n-2}\lra X_{n-3}\lra \cdots\lra X_1\lra X_0\lra {}_AN\otimes_BV\lra 0$$with $X_j\in\add(_AX)$ such that $\Hom_A(X',-)$ preserves exactness of $(*)$ for $X'\in \add(_AX)$. Since $M_A$ is projective, we get an exact sequence
$$ (**)\quad 0\lra {}_BM\otimes_AX_{n-2}\lra {}_BM\otimes_AX_{n-3}\lra \cdots\lra {}_BM\otimes_AX_1\lra {}_BM\otimes_AX_0\lra {}_BM\otimes_AN\otimes_BV\lra 0.$$ Clearly  $_BM\otimes_AX_j\in\add(_BM\otimes_AX)$. To see that $\Hom_B(_BM\otimes_AX, -)$ preserves exactness of $(**)$, we note the isomorphisms
$$ \Hom_B(M\otimes_AX, M\otimes_AX_j)\simeq \Hom_A(M\otimes_AX, \Hom_A(_AN_B,X_j))\simeq \Hom_A(_AN\otimes_BM\otimes_AX,X_j).$$ Note that $_AN\otimes_BM_A\simeq \End_B(M)$ as $A$-$A$-bimodules. Since $_AN\otimes_BM_A$ is centrally projective over $A$, that is, ${}_AN\otimes_BM_A\in \add(_AA_A)$, we have $$_AN\otimes_BM_A\otimes_AX_j\in \add(_AA\otimes_AX_j)\in \add(_AX).$$
Hence, in the commutative diagram
$${\footnotesize \xymatrix{
0\ar[r] &(M\otimes_AX, M\otimes_AX_{n-2})\ar[r]\ar[d]^-{\simeq}& \cdots \ar[r]&(M\otimes_AX,M\otimes_AX_0)\ar[d]^-{\simeq}\ar[r]& (M\otimes_AX,M\otimes_AN\otimes_BV)\ar[d]^-{\simeq}\ar[r]& 0\\
0\ar[r] &(N\otimes_BM\otimes_BX,X_{n-2})\ar[r] &\cdots\ar[r]&(N\otimes_BM\otimes_BX,X_0)\ar[r] &(N\otimes_BM\otimes_BX,N\otimes_BV)\ar[r]&0,
}}$$
the bottom row is exact. This means that the top row is also exact and that the projective dimension of the $\End_B(M\otimes_AX)$-module $\Hom_B(M\otimes_AX, M\otimes_AN\otimes_AV)$ is at most $n-2$. Since $B$ is $M$-separable over $A$, that is, $_BB_B\in \add(_B(M\otimes_AN)_B)$, we have $_BV\simeq {}_BB\otimes_BV\in \add\big(_B(M\otimes_AN)\otimes_BV\big)$. Thus the $\End_B(M\otimes_AX)$-module $\Hom_B(M\otimes_AX, V)$, as a direct summand of the module $\Hom_B(M\otimes_AX, M\otimes_AN\otimes_AV)$, has projective dimension at most $n-2.$ Hence $\gd(\End(_BM\otimes_AX))\le n$ by Lemma \ref{lem1}. This implies
$$\rd(B)\le \gd(\End_B(M\otimes_AX))\le n=\rd(A). \; \; \square $$

\medskip
{\bf Proof of Theorem \ref{sym-equ}}. Suppose that Artin algebras $A$ and $B$ are symmetric separably equivalent through  bimodules $_AP_B$ and $_BQ_A$. Then, by \cite[Theorem 4.1]{Kadison2019}, there is a split, separable, Frobenius extension $A\subseteq \End_{B^{op}}(P_B)$, with $P_B$  a projective generator for $B^{op}\modcat$. Dually, there is a split, separable, Frobenius extension $B\subseteq \End_{A^{op}}(Q_A)$, with $Q_A$  a projective generator for $A^{op}\modcat$. Generally, for a semisimple extension $\Gamma\subseteq \Lambda$ of Artin algebras, if $\Gamma$ is semisimple, then $\Lambda$ is semisimple. Thus $A$ is semisimple if and only if $B$ is semisimple. So we may assume that both $A$ and $B$ are non-semisimple.

We shall prove that if the $A$-$A$-bimodule $P\otimes_BQ$ is centrally projective over $A$, then $\rd(A)\le \rd(B)$. In fact, we verify that $_AP_B$ satisfies the conditions in Theorem \ref{frob-bim}. First, according to \cite[Theorem 4.1]{Kadison2019}, $_AP_B$ and $_BQ_A$ are Frobenius bimodules. Second, $_AP$ is a projective generators. This follows from the separability: $_AA\in \add(_AP\otimes_BQ)\subseteq \add(_AP\otimes_BB)=\add(_AP)$. Clearly, $A$ is $P$-separable over $A$. Since $\End_{B^{op}}(P_B)\simeq P\otimes_BQ$ as $A$-$A$-bimodules, $_A\End_{B^{op}}(P_B)_A$ is centrally projective over $A$ by assumption on $P\otimes_BQ$. Thus $\rd(A)\le \rd(B)$ by Theorem \ref{frob-bim}.

Dually, the bimodule $_BQ_A$ fulfills the conditions of Theorem \ref{frob-bim}, we get $\rd(B)\le \rd(A)$. Thus $\rd(A)=\rd(B)$. $\square$

\begin{Rem}{\rm (1) Separable equivalent algebras have the same finitistic dimensions, while symmetric separable equivalent algebras have the same dominant dimensions. Recall that the finitistic dimension of an Artin algebra $A$, denoted by $\fd(A)$, is the supremum of the projective dimensions of finitely generated modules of finite projective dimensions. In fact, suppose that bimodules $_AP_B$ and $_BQ_A$ define such an equivalence. Then, for an $A$-module $X$, we have $\pd(_AX)\le \pd(_AP\otimes_BQ\otimes_AX)\le\pd(_BQ\otimes_AX)\le \pd(_AX)$. This implies $$\fd(A) = \mbox{sup}\{\pd(_BQ\otimes_AX)\mid \pd(_AX)<\infty\}\le \fd(B).$$ Dually, $\fd(B)\le \fd(A).$ Thus $\fd(A)=\fd(B)$. This argument also shows that global dimensions are invariant under separable equivalences, as pointed in \cite{kadison95}. For dominant dimensions, we note that the functors of the adjoint triple $(P\otimes_B-,Q\otimes_A-, P\otimes_B-)$ preserve projective-injective modules. This implies the equality of dominant dimensions of $A$ and $B$ if they are symmetric separably equivalent. For the precise definition of dominant dimensions, we refer to \cite{x1} and the references therein.

(2) If two Artin algebras $A$ and $B$ are symmetric separably equivalent, then one of them is a Gorenstein algebra of dimension $d$ (that is, $\id(_AA)= d =\id(A_A)$) if and only if so is the other. Indeed, the adjoint pairs $(P\otimes_B-,Q\otimes_A-)$ and $(Q\otimes_A-,P\otimes_B-)$ yield that $_AP\otimes_BY$ and $_BQ\otimes_AX$ are injective modules if $_AX$ and $_BY$ are injective modules because, for example, in case of $_B Y$ being injective, the functor
$$\Hom_A(-, P\otimes_BY)\simeq \Hom_A(-, \Hom_B(Q,Y))\simeq \Hom_B(Q\otimes_A-, Y)=\Hom_B(?,Y)\; (Q\otimes_A-)$$is exact. As in (1), we get $\id(_BY)=\id(_AP\otimes_BY)$ and $\id(_AX)=\id(_BQ\otimes_AX)$. This gives the desired conclusion. Consequently, for $d=0$, we see that $A$ is a self-injective algebra if and only if so is $B$.
}\end{Rem}

\medskip
Now, we apply the theorems to Frobenius extensions and get results of two types.

\begin{Koro} \label{frob-ext}
$(1)$ If $B\subseteq A$ is a separable extension of Artin algebras such that $_BA_B$ is centrally projective over $B$, then $\rd(A)\le \rd(B)$. Moreover, if $C$ is an intermediate algebra between $B$ and $A$ such that $_CA_C$ is centrally projective over $C$, then $\rd(C)\le \rd(B).$

$(2)$ If $B\subseteq A$ is a split, $H$-separable extension of Artin algebras, that is, $A\otimes_BA\in \add(_AA_A)$, and if $_BA$ is projective, then $\rd(B)=\rd(\End(_BA))\le \rd(A).$
\end{Koro}

{\it Proof.} (1) Since the extension $B\subseteq A$ is separable, it is also a semisimple extension. Thus, if $B$ is semisimple, then so is $A$ because semisimple extensions preserve semisimplicity, and $\rd(B)=\rd(A)$. So we assume that $B$ is non-semisimple, that is, $\rd(B)$ is at least $2$. We consider the bimodule ${}_AA_B$.  Clearly, $^*(_AA_B)=\Hom_A(_AA_B,A)={}_BA_A$ and $_AA$ is a generator for $A$-mod. Since the extension is separable, we have $_AA_A\in \add(_AA\otimes_BA_A)$, that is, $A$ is $_AA_B$-separable over $B$. By assumption, $_BA_B\simeq {}_BA\otimes_AA_B$ is centrally projective over $B$.  Then, by Lemma \ref{sep-over}(3) (or \cite[Theorem 2]{sugano1}), $B\subseteq A$ is a Frobenius extension, and therefore the $A$-$B$-bimodule $_AA_B$ is a Frobenius bimodule. Hence the first statement in (1) follows from Theorem \ref{frob-bim}.

To get the last statement in (1), we note that, if $A$ is centrally projective over $B$, then it is shown in \cite[Proposition 13, p.206]{sugano1} and the  proof there that $B\subseteq C$ is a separable extension and
$_BC_B$ is centrally projective over $B$. Thus $\rd(C)\le\rd(B)$.

(2) Since $_BA$ is projective and the extension is split, it follows from Lemma \ref{sep-over}(2)(a) or \cite[Theorem 7(1)]{sugano1} that $A\subseteq \End(_BA)$ is a separable extension. As $B\subseteq A$ is an $H$-separable extension, that is, $_AA\otimes_BA_A\in \add(_AA_A)$, we obtain the $A$-$A$-bimodule isomorphism and relations
$$_A\End(_BA)_A = \Hom_B(_BA_A,{}_BA_A)\simeq \Hom_A(_AA\otimes_BA, {}_AA)\in\add(_A\Hom_A(_AA, {}_AA)_A)\subseteq \add(_AA_A), $$
that is, the $A$-$A$-bimodule $\End(_BA)$ is centrally projective over $A$. Thus the extension $A\subseteq \End(_BA))$ satisfies all conditions of Corollary \ref{frob-ext}(1) and therefore $\rd(\End(_BA))\le \rd(A).$  Note that the splitting condition shows that $_BB\in \add(_BA)$ and $_BA$ is a projective generator for $B$-mod. Hence $B$ and $\End(_BA)$ is Morita equivalent. This implies that $\rd(B)=\rd(\End(_BA))\le \rd(A).$
$\square$

Remark that if $A_B$ is also projective in Corollary \ref{frob-ext}(2), then, by Lemma \ref{sep-over}(2), $ B\subseteq A$ is a Frobenius extension. Thus $_BA_A$ is a Frobenius bimodule and $_BA$ is a generator for $B\modcat$, and $B$ is $_BA_A$-separable over $A$. This can be seen by choosing a Frobenius system $(E, x_i, y_i)$ and checking that $E$ is split surjective as $B$-bimodules (see \cite[p.14]{kadison}). By the assumption of $H$-separability, $A\otimes_BA$ is centrally projective over $A$. Thus all of the conditions of Theorem \ref{frob-bim} are met. Whence (2) follows.

\medskip
There are examples of extensions satisfying all conditions of Corollary \ref{frob-ext}(2). For example, for any Artin algebra $B$ with an automorphism $\rho$, let $B[X,\rho]$ be the skew polynomial algebra with $Xb=(b)\rho\, X$. If $f\in B[X,\rho]$ is a monic polynomial of degree $n$ such that $fB[X,\rho]=B[X,\rho]f$, then there is a skew polynomial algebra $B[x,\rho]$, defined as the quotient algebra $B[X,\rho]/fB[X,\rho]$, of degree $n$ with a free $B$-basis $\{1, x, \cdots, x^{n-1}\}$, where $x=X+ fB[X,\rho]$. Let $C$ be the center of $B$, and $B^{\rho}$ be the $\rho$-fixed subalgebra of $B$. If $f=x^n-u$ for an invertible element $u\in B^{\rho}$, then it was shown in \cite{sx} that $B[x,\rho]$ is an $H$-separable extension of $B$ if and only if $C$ is a $\rho$-Galois extension of $C^{\rho}$. Clearly, the extension $B\subseteq B[x,\rho]$ is split. Recall that $B$ is called a \emph{$\rho$-Galois extension} of $B^{\rho}$ if there exist elements $\{x_i, y_i\mid 1\le i\le m\}$ in $B$ for some positive integer $m$ such that $\sum_i x_iy_i=1$ and $\sum_ix_i(y_i)\rho^j=0$ for $1\le j\le m.$

\begin{Koro} \label{H-sep-cen-pro} Let $B\subseteq A$ be an $H$-separable extension of Artin algebras such that the $B$-$B$-bimodule $_BA_B$ is centrally projective over $B$. Then, for any $B$-module $_BY\in B\modcat$, we have $\rd(\End_A(A\otimes_BY))\le \rd(\End_B(Y))$.% if $\End_B(Y)$ is non-semisimple.
\end{Koro}

{\it Proof.} We cite a result in \cite{nak}: If $S\subseteq R$ is an $H$-separable extension of rings such that $_SR_S$ is centrally projective over $S$, then so is the extension $\End_S(V)\subseteq \End_R(R\otimes_SV)$, $f\mapsto id_A\otimes f$, for any $S$-module $V$. Thus Corollary \ref{H-sep-cen-pro} follows immediately from
Corollary \ref{frob-ext}(1) since $H$-seprable extensions are separable extensions.

\begin{Koro} \label{sep-ten} Let $k$ be a field. If $B\subseteq A$ is a separable extension of finite-dimensional $k$-algebras, such that $_BA_B$ is centrally projective over $B$, then $\rd(A\otimes_kC) \le \rd(B\otimes_kC)$ for all finite-dimensional $k$-algebra $C.$
\end{Koro}

{\it Proof.} By \cite[Corollary 2.8]{hirata-sugano}, as $B\subseteq A$ is a separable extension, $B\otimes_kC\subseteq A\otimes_kC$ is a separable extension for any $k$-algebra $C$. Since $_BA_B$ is centrally projective over $B$, ${}_{B\otimes_kC}(A\otimes_kC)_{B\otimes_kC}\simeq {}_{B\otimes_k C}\big({}_BA_B\otimes_k {}_CC_C\big)_{B\otimes_kC}\in \add(_BB_B\otimes_k{}_CC_C) = \add(_{B\otimes_kC}(B\otimes_kC)_{B\otimes_kC}).$
That is, $_{B\otimes_kC}(A\otimes_kC)_{B\otimes_kC}$ is centrally projective over $B\otimes_kC$. Thus
Corollary \ref{sep-ten} follows from Corollary \ref{frob-ext}(1). $\square$

\medskip
Remark that we always have $\rd(A\otimes_kB)\le \rd(A)+\rd(B)$ for finite-dimensional algebras $A$ and $B$ over a perfect field $k$ (see \cite{x2}).

As a consequence of Corollary \ref{frob-ext}, we have the following.

\begin{Koro} Let $B\subseteq A$ be a split, separable£¬ Frobenius extension of Artin algebras. If $_BA_B\in \add(_BB_B)$ and $_A(A\otimes_BA)_A\in \add(_AA_A)$, then $\rd(A)=\rd(B)$. \end{Koro}

\medskip
For a finite-dimensional Hopf-algebra $H$ over a field $k$,
a left $H$-module $k$-algebra $A$, and an invertible normal two cocycle $\sigma: H\otimes_kH\ra A$, there is defined a crossed product $A\#_{\sigma}H$ of $A$ and $H$ with respect to $\sigma$, which is an associative $k$-algebra with the underlying $k$-space $A\otimes_kH$ and the identity $1\# 1$. For the precise definition of crossed products, one may refer to, for example, the book of Montgomery \cite[\S 7, p.101]{montgomery} or the original reference therein. For representation dimensions of crossed products, we obtain the following result.

\begin{Koro}\label{hopf-algebra} Let $H$ be a finite-dimensional, cocommutative Hopf-algebra over a field, $A$
a left $H$-module algebra and $\sigma$ an invertible normal two cocycle. If $A$ has a central element of trace $1$ , then $\rd(A\#_{\sigma}H) \le \rd(A)$.
In particular, $\rd(A\# H) \le \rd(A)$, where $A\# H$ is the smash product of $A$ with $H$.
\end{Koro}

{\it Proof.} Let us consider the canonical extension $A\ra A\#_{\sigma}H, a\mapsto a\#1$.
By \cite[Corollary 7.2.11, p.111]{montgomery}, $A\#_{\sigma}H$ is a direct sum of $\dim_k(H)$ copies of $A$ as right and left $A$-modules, thus as $A$-$A$-bimodules. This means that the $A$-$A$-bimodule $_A(A\#_{\sigma}H)_A$ is central over $A$. Note that a finite-dimensional Hopf algebra over a field has antipode always invertible. Since $H$ is cocommutative and $A$ has a central trace $1$ element, the extension $A\subseteq$ $A \#_{\sigma} H$ is separable by \cite[Theorem 1.11]{cf}. Hence Corollary \ref{hopf-algebra} follows from Corollary \ref{frob-ext}(1). The last statement follows from the first one.
$\square$

\medskip
{\bf Proof of Corollary \ref{int-cor}}. (2) is shown in Corollary \ref{hopf-algebra}. Note that for self-injective Artin algebras $A$, there holds always $\rd(A)\le LL(A)$, where $LL(A)$ denotes the Loewy length of $A$. Further, if $A$ is a finite-dimensional algebra over a field $k$, then $LL(A) \le \dim_k(A)$.

It remains to show (1).  First, we suppose that $H$ is a normal subgroup of $G$.

In this case, by the Mackey formula or double coset decomposition of $G$ with respect to $H$, $_{k[H]}k[G]_{k[H]}\in \add(_{k[H]}k[H]_{k[H]})$ by \cite[p. 61]{benson}.
For an Artin algebra $k$ and any subgroup $H$ of a finite group $G$ such that $[G:H]$ is invertible in $k$, it is known that $k[H]\subseteq k[G]$ is a separable, Frobenius extension with the Frobenius homomorphism $$E: k[G]\lra k[H], \quad \sum_{g\in G}\lambda_gg\mapsto \sum_{g\in H}\lambda_gg$$(for example, see \cite[Example 3.1]{kadison95} for details). Since $_{k[H]}k[G]_{k[H]}\in \add(_{k[H]}k[H]_{k[H]})$, the $k[H]$-$k[H]$-bimodule $k[G]$ is centrally projective over $k[H]$,  and $\rd(k[G])\le \rd(k[H])$ by Corollary \ref{frob-ext}(1).  This proves (1) for the case of $H$ being normal subgroup of $G$.

Next, we consider any subgroup $H$ of $G$ such that $[G:H]$ is invertible in $k$. If $H$ is a normal subgroup of $G$, then the conclusion follows. Suppose that $H$ is not normal in $G$ and we then define $H_1:=N_G(H)$, the normalizer of $H$ in $G$. Then $H\lneq H_1$ and $[H_1:H]$ is invertible in $k$. Thus $\rd(H_1)\le \rd(H)$, as just proved. Note that $[G:H_1]$ is invertible in $k$. If $H_1$ is normal in $G$, then $\rd(k[G]) \le \rd(k[H_1])\le \rd(k[H])$, as desired. If $H_1$ is not normal in $G$, then we define $H_2:=N_G(H_1)$. Thus $H\lneq H_1\lneq H_2$ and $\rd(k[H_2])\le \rd(k[H_1]).$ If $H_2$ is normal in $G$, then $\rd(k[G])$ $ \le \rd(k[H_2]) \le \rd(k[H_1] \le \rd(k[H])$. If $H_2$ is not normal in $G$, we repeat the foregoing argument. Since $G$ is a finite group, there exists some $H_i$ which is the normalizer of $H_{i-1}$ in $G$ such that $H_i$ is a normal subgroup of $G$. Since $H_i$ contains $H$, the index $[G:H_i]$ is invertible in $k$. Hence $\rd(k[G])\le \rd(k[H_i])\le \cdots \le\rd(H_1)\le\rd(k[H]).$ Thus the proof of (1) is completed.

Now, let $k$ be a field and $H$ is a subgroup of $G$ with $[G:H]$ invertible in $k$. We assume that $H$ is normal in $G$ and write $k[G]$ as $k[H]\#_{\sigma}k[G/H] $ and use the Blattner-Montgomery duality (see \cite[Corollary 9.4.17]{montgomery}) as did in \cite{sun}: $(A\#_{\sigma}H)\#H^*\simeq M_n(A)$, where $H^*$ is the $k$-dual Hopf algebra and $n=\dim_k(H).$ Thus it follows from (2) that
$$\begin{array}{rl} \rd(k[H])=\rd(M_n(k[H])& =\rd\big((k[H]\#_{\sigma}k[G/H])\#(G/H)^*\big)\\ & \le \rd\big(k[H]\#_{\sigma}(G/H)\big) \\ & \le \rd(k[H]),\end{array} $$ That is, $\rd(k[G])=\rd(k[H]).$

If $H$ is not normal in $G$, we may use the argument in the proof of (1). For instance, we define $H_1:=N_G(H)$, the normalizer of $H$ in $G$. Then $H\lneq H_1$ and $[H_1:H]$ is invertible in $k$. Thus $\rd(H_1) = \rd(H)$, as just proved. Continuing this procedure, we get some $H_i$ which is the normalizer of $H_{i-1}$ in $G$ such that $H_i$ is a normal subgroup of $G$. Since $H_i$ contains $H$, the index $[G:H_i]$ is invertible in $k$. Hence $\rd(k[G]) = \rd(k[H_i])= \cdots = \rd(H_1)=\rd(k[H]).$ (3) is proved.
$\square$

Remark that, if $H$ is a normal subgroup of a finite group $G$ such that $[G:H]$ is invertible in a field $k$, then the group algebra $k[G/H]$ is semisimple by Maschke's theorem and cosemisimple (see \cite[2.4.1, p.25]{montgomery}). Thus, in this case, the equality of $\rd(k[G])=\rd(k[H])$ follows also from \cite[Corolary 3.5]{sun}.

%Second proof Corollary \ref{int-cor}(1): the restriction of $E$ to $k[H]$ is the identity map on $k[H]$. Therefore $k[H]$ and $k[G]$  are symmetric separably equivalent Artin algebras. This is even true for $k$ to be any ring and for $G$ to be any group such that $[G:H]$ is finite and invertible in $k$ (see \cite[p.348]{Kadison2019}).
%Since $k[H]$-$k[H]$-bimodule $k[G]$ is centrally projective over $k[H]$,Corollary \ref{int-cor}(1) follows from the proof of Theorem \ref{sym-equ} immediately.

\begin{Koro} If $k$ is a field of positive characteristic $p$ and $P$ is a Sylow $p$-subgroup of a finite group $G$, then $\rd(k[G]) = \rd(k[P])\le |P|.$ \label{p-gr}
\end{Koro}

Corollary \ref{p-gr} shows that representation dimensions of group algebras over a field for arbitrary finite groups are reduced to the ones for finite $p$-groups. For example, we have the following consequence.

\begin{Koro} \label{iso-syl} Let $k$ be a perfect field of positive characteristic $p$, and let $G$ be finite groups.
If a Sylow $p$-subgroup of $G$ is of the form $C_{p^{n_1}}\times\cdots\times C_{p^{n_s}}$, where $n_j\ge 1$ and $C_n$ is the cyclic group of order $n$, then $\rd(k[G])\le 2s.$
\end{Koro}

{\it Proof.} By Corollary \ref{p-gr}, $\rd(k[G])=\rd(k[P])=\rd(k[C_{p^{n_1}}\times\cdots\times C_{p^{n_s}}])$. While $k[C_{p^{n_1}}\times\cdots\times C_{p^{n_s}}]\simeq k[C_{p^{n_1}}]\otimes_k\cdots\otimes_kk[C_{p^{n_s}}]$ and $\rd(k[C_{p^n}])=2$ for $n\ge 1$, it then follows from \cite[Theorem 3.5]{x2} that
$\rd(k[G])\le 2s.$ $\square$

\medskip
Now we pass to considering Frobenius extensions of quotient algebras. First, we mention the following lemma which is ready to prove.

\begin{Lem}\label{tensorchange} Let $I$ be an ideal in a ring $R$ and $\bar{R}:=R/I$. Let $X$ and $U$ be right $R$-modules and $UI=0$. If ${}_RV$ is an $R$-module with $IV=0$, then $U\otimes_RV\simeq U\otimes_{\bar{R}}V$ and $X\otimes_RV\simeq (X/XI)\otimes_RV\simeq  (X/XI)\otimes_{\bar{R}}V.$
\end{Lem}

\begin{Koro} Let $B\subseteq A$ be a separable extension of Artin algebras such that $_BA_B$ is centrally projective over $B$. If $I$ is an ideal in $A$ such that $A(I\cap B)+(I\cap B)A=I$, then $\rd(A/I)\le \rd(B/(I\cap B))$.
\label{factor}
\end{Koro}

{\bf Proof.} Set $\bar{A}:=A/I$ and $\bar{B}:=B/(I\cap B)$. We will show that $\bar{B}\subseteq \bar{A}$ satisfies all conditions in Corollary \ref{frob-ext}(1). It is known (for example, see \cite{hirata-sugano}) that an extension $S\subseteq R$ of rings is separable if and only if there are elements $x_i, y_i\in R$, $1\le i\le m$, such that $\sum_i x_iy_i=1$ and $\sum_i a x_i\otimes y_i=\sum_i x_i\otimes y_ia$ in $R\otimes_S R$ for all $a\in R.$ Thus, for the separable extension $B\subseteq A$, there are $2n$ elements $a_i,b_i\in A$ such that $\sum_i a_ib_i=1$ and $\sum_i a a_i\otimes b_i=\sum_i a_i\otimes b_ia.$ Let $^{-}: A\ra \bar{A}, a\mapsto \bar{a}=a+I$ be the canonical projection. Then $\sum_i \bar{x}_i\bar{y}_i=1.$ By Lemma \ref{tensorchange}, under the map $^-\otimes_B^-$, we also have $\sum_i \bar{a} \bar{x}_i\otimes_{\bar{B}} \bar{y}_i=\sum_i \bar{x}_i\otimes_{\bar{B}} \bar{y}_i\bar{a}$ in $\bar{A}\otimes_{\bar{B}}\bar{A}$ for all $a\in A$. Thus the extension $\bar{B}\subseteq \bar{A}$ is separable.

To see that the $\bar{B}$-$\bar{B}$-bimodule $\bar{A}$ is centrally projective over $\bar{B}$, we consider  $B$-$B$-bimodules as $B^e$-module, where $B^e:=B\otimes_kB^{op}$ is the enveloping algebra of $B$. Then  $\bar{B}^e= B^e/J$ with $J= (I\cap B)\otimes_kB^{op}+B\otimes_k(I^{op}\cap B^{op})$.  Assume $_BA_B\oplus {}_BY_B\simeq \bigoplus_{j=1}^m{}_BB_B$. Then, applying the functor $(B^e/J)\otimes_B-$ to this isomorphism, we get
$$ {}_{\bar{B}}\big(A/(A(I\cap B) +(I\cap B)A)\big)_{\bar{B}}\oplus (B^e/J)\otimes_{B^e} Y\simeq \bigoplus_{j=1}^m{}_{\bar{B}}\big(B/(I\cap B)\big)_{\bar{B}}$$
Since $I=A(I\cap B)+(I\cap B)A$, we obtain $_{\bar{B}}\bar{A}_{\bar{B}}\in \add(_{\bar B}\bar{B}_{\bar{B}}).$ Thus Corollary \ref{factor} follows from Corollary \ref{frob-ext}(1). $\square$

Let us end this section by a couple of  examples.

(1) Suppose $k$ is a field of characteristic $2$. Let $A_4$ and $V_4$
be the alternative group of degree $4$ and Klein $4$-group, respectively. Then $k[A_4]$ and $k[V_4]$ are not stably equivalent (see \cite{peacock}), but they are symmetric separably equivalent (see the proof of Corollary \ref{int-cor}). Clearly, $V_4$ is a normal Sylow $2$-subgroup of $A_4$. By Corollary \ref{int-cor}(1), $\rd(k[A_4])=\rd(k[V_4]).$ Note that both $k[A_4]$ and $k[V_4]$ is representation-infinite by a result of Higman, which says that the group algebra $k[G]$ of a finite group $G$ over a field $k$ of characteristic $p>0$ is representation-finite if and only if Sylow $p$-groups of $G$ are cyclic. Since $V_4\simeq C_2\times C_2$ and $k[V_4]$ is of Loewy length $3$, we have $\rd(k[A_4])=\rd(k[V_4])=3$.

(2) Let $p\ne q$ be two prime numbers, and let $S_q$ and $A_q$ be the symmetric and alternative groups of degree $q$, respectively. If $k$ is a field of characteristic $p$, then $\rd(k[S_q])=\rd(k[S_{q-1}])$ and $\rd(k[A_q])=\rd(k[A_{q-1}]).$ This is a consequence of Corollary \ref{p-gr}.

\medskip
{\bf Acknowledgements.} The author is greatly indebted to Professor Shenglin Zhu at Fudan University for offering his help with the literature on Hopf algebras, and for discussing separable extensions. The research work of the author was partially supported by NSFC (12031003) and BJNSF (1192004).

\medskip
{\footnotesize

Changchang Xi, School of Mathematical Sciences, Capital Normal University, 100048 Beijing, China; and School of Mathematics and Statistics, Central China Normal University, 430079 Wuhan, China

{\tt Email: xicc@cnu.edu.cn}

}

\end{document}